\documentclass[12pt]{article}
\usepackage{amsfonts,amsmath,amssymb}
\usepackage{enumerate}
\usepackage{theorem}
\numberwithin{equation}{section}

\newtheorem{thm}{Theorem}[section]
\newtheorem{prop}{Proposition}[section] 
\newtheorem{lem}{Lemma}[section] 
\newtheorem{cor}{Corollary}[section]
\newtheorem{Def}{Definition}[section]
\newtheorem{rem}{Remark}[section]

\newenvironment*{Proof}{{\bf Proof.}}

\newcommand{\fa}{\forall}
\newcommand{\dis}{\displaystyle}

\newcommand{\cm}{\mathcal{M}}
\newcommand{\ce}{\mathcal{E}}
\newcommand{\ct}{\mathcal{T}}
\newcommand{\cf}{\mathcal{F}}
\newcommand{\ld}{\ldots}

\newcommand{\ra}{\rightarrow}

\newcommand{\al}{\alpha}
\newcommand{\bi}{\beta}

\newcommand{\Ga} {{\varGamma}}

\newcommand{\de}{\delta }

\newcommand{\De} {{\varDelta}}
\newcommand{\e}{\varepsilon }
\newcommand{\f}{\phi}

\newcommand{\vthi}{\vartheta }

\newcommand{\la}{\lambda }
\newcommand{\mi}{\mu }

\newcommand{\sm}{\smallsetminus}

\newcommand{\oo}{\omega}

\newcommand{\R}{\mathbb{R}}

\newcommand{\Z}{\mathbb{Z}}
\newcommand{\N}{\mathbb{N}}

\newcommand{\ssum}{\sum\limits}
\newcommand{\qs}{$\quad\square$}

\begin{document}
\title{\bf Optimal weak type estimates for dyadic-like maximal operators}
\author{Eleftherios N. Nikolidakis}
\date{}
\maketitle

\begin{abstract}
We provide sharp weak estimates for the distribution function of
$\cm\f$ when on $\f$ we impose $L^1$, $L^q$ and $L^{p,\infty}$
restrictions. Here $\cm$ is the dyadic maximal operator associated
to a tree $\ct$ on a non-atomic probability measure space.
\end{abstract}
\section{Introduction} 
\noindent

The dyadic maximal operator on $\R^n$ is defined by
\begin{eqnarray}
\cm_d\f(x)=\sup\bigg\{\frac{1}{|Q|}\int_Q|\f(u)|du:\; x\in
Q,\;Q\subseteq\R^n \ \ \text{is a dyadic cube}\bigg\}
\label{eq1.1}
\end{eqnarray}
for every $\f\in L^1_{loc}(\R^n)$ where the dyadic cubes are those
formed by the grids $2^{-N}\Z^n$ for $N=1,2,\ld\;.$

It is well known that it satisfies the following weak type
(\ref{eq1.1}) inequality
\begin{eqnarray}
|\{x\in\R^n:\cm_d\f(x)>\la\}|\le\frac{1}{\la}\int_{\{\cm_d\f>\la\}}|\f(u)|du
\label{eq1.2}
\end{eqnarray}
for every $\f\in L^1(\R^n)$ and every $\la>0$.

Using (\ref{eq1.1}) we easily get the following $L^p$ inequality
\begin{eqnarray}
\|\cm_d\f\|_p\le\frac{p}{p-1}\|\f\|_p  \label{eq1.3}
\end{eqnarray}
for every $p>1$ and every $\f\in L^p(\R^n)$, which is proved to be
best possible (see \cite{2}, \cite{3} for the general martingales
and \cite{10} for the dyadic ones).

A way of studying the dyadic maximal operator is the introduction
of the so called Bellman functions (see \cite{8}).

Actually, we define for every $p>1$
\begin{eqnarray}
B_p(f,F)=\sup\bigg\{\frac{1}{|Q|}\int_Q(\cm_d\f)^p:Av_Q(\f^p)=F,
\ \ Av_Q(\f)=f\bigg\}  \label{eq1.4}
\end{eqnarray}
where $Q$ is a fixed dyadic cube, $\f$ is nonnegative in $L^p(Q)$
and $f,F$ are such that $0<f^p\le F$.

$B_p(f,F)$ has been computed in \cite{5}. In fact it has been shown
that $B_p(f,F)=F\oo_p(f^p/F)^p$ where
$\oo_p:[0,1]\ra\Big[1,\frac{p}{p-1}\Big]$ is the inverse function of
\[
H_p(z)=-(p-1)z^p+pz^{p-1}.
\]
Actually this has been proved in a much more general setting of
tree like maximal operators on non-atomic probability spaces. The
result turns out to be independent of the choice of the measure
space.

The study of these operators has been continued in \cite{7} where
the Bellman functions of them in the case $p<1$ have been computed.

Actually, as in \cite{5} and \cite{7} we will take the more
general approach. So for a tree $\ct$ on a non atomic probability
measure space $X$, we define the associated dyadic maximal
operator, namely
\[
\cm_\ct\f(x)=\sup\bigg\{\frac{1}{\mi(I)}\int_I|\f|d\mi:\;x\in
I\in\ct\bigg\}
\]
for every $\f\in L^1(X,\mi)$.

It is now known that $\cm_\ct:L^{p,\infty}\ra L^{p,\infty}$ is a
bounded operator satisfying
\begin{eqnarray}
\|\cm_\ct\f\|_{p,\infty}\le|||\f|||_{p,\infty}.  \label{eq1.5}
\end{eqnarray}

It is now interesting to see what happens if we replace the
$L^p$-norm of $\f$ in (\ref{eq1.4}) by it's $L^{p,\infty}$-norm,
$|||\cdot|||_{p,\infty}$, given by
\begin{align*}
|||\f|||_{p,\infty}=\sup\Big\{\mi(E)^{-1+\frac{1}{p}}\int_E|\f|d\mi:
\;&E\ \ \text{measurable subset of $X$ such that}\\
&\mi(E)>0\bigg\}.
\end{align*}
It is well known that $|||\cdot|||_{p,\infty}$ is a norm on
$L^{p,\infty}$ equivalent to the quasi norm $\|\cdot\|_{p,\infty}$
defined by
\[
\|\f\|_{p,\infty}=\sup\Big\{\la\mi(\{\f\ge\la\})^{1/p}: \;
\la>0\Big\}.
\]
In fact
\[
\|\f\|_{p,\infty}\le|||\f|||_{p,\infty}\le\frac{p}{p-1}\|\f\|_{p,\infty},
\ \ \fa\;\f\in L^{p,\infty}
\]
as can been seen in \cite{4}.

In fact in \cite{9} it is proved that (\ref{eq1.5}) is sharp
allowing every value for the $L^1$-norm of $\f$.

In the present paper we compute the following function
\begin{align}
S(f,A,F,\la)=\sup\bigg\{&\mi(\{\cm_\ct\f\ge\la\}):\;\f\ge0,\;
\int_X\f d\mi=f, \nonumber \\
&\int_X\f^qd\mi=A, \ \ |||\f|||_{p,\infty}=F\bigg\}
\label{eq1.6}
\end{align}
for every $\la>0$, $(f,A,F)$ on the domain of the extremal problem
and $q$ fixed such that $1<q<p$. That is we provide improvements of
(\ref{eq1.3}) given additionally $L^q$ and $L^{p,\infty}$
restrictions on $\f$.

From this we obtain as a corollary that
\begin{equation}
\sup\bigg\{\|\cm_\ct\f\|_{p,\infty}:\;\f\ge0,\;\int_X\f d\mi=f,\;
\int_X\f^qd\mi=A, \ \ |||\f|||_{p,\infty}=F\bigg\}=F
\label{eq1.7}
\end{equation}
that is (\ref{eq1.5}) is sharp allowing every value of the
integral and the $L^q$-norm of $\f$, for a fixed $q$ such that
$1<q<p$. As a matter of fact we prove that the supremum in both
cases (\ref{eq1.6}) and (\ref{eq1.7}) is attained. These estimates
are provided in Section 4, while in Section 3 the domain of the
extremal problem is found. On Section 2 we give some preliminaries
needed during this paper.

Finally we mention that all the above estimates are independent of
the measure space and the tree $\ct$.
\section{Preliminaries} 
\noindent

Let $(X,\mi)$ be a non-atomic probability measure space. We state
the following lemma which can be found in \cite{1}.
\begin{lem}\label{lem2.1}
Let $\f:(X,\mi)\ra\R^+$ and $\f^\ast$ the decreasing rearrangement
of $\f$, defined on $[0,1]$. Then
\[
\int^t_0\f^\ast(u)du=\sup\bigg\{\int_E\f d\mi:\; E \;\text{measurable subset of}\;\; X \;
\text{with}\; \mi(E)=t \bigg\}
\]
for every $t\in[0,1]$, with the supremum in fact attained.  \qs
\end{lem}

We prove now the following:
\begin{lem}\label{lem2.2}
Let $\f:X\ra\R^+$ be measurable and $I\subseteq X$ be measurable
with $\mi(I)>0$. Suppose that $\frac{1}{\mi(I)}\int\limits_I\f
d\mi=s$. Then for every $t$ such that $0<t\le\mi(I)$ then exists a
measurable set $E_t\subseteq I$ with $\mi(E_t)=t$ and
$\frac{1}{\mi(E_t)}\int\limits_{E_t}\f d\mi=s$.
\end{lem}
\begin{Proof} Consider the measure space $(I,\mi/I)$ and let
$\psi:I\ra\R^+$ be the restriction of $\f$ on $I$ that is
$\psi=\f/I$. Then if $\psi^\ast:[0,\mi(I)]\ra\R^+$ is the decreasing
rearrangement of $\psi$, we have that
\begin{eqnarray}
\frac{1}{t}\int^t_0\psi^\ast(u)du\ge\frac{1}{\mi(I)}\int^{\mi(I)}_0\psi^\ast
(u)du=s\ge\frac{1}{t}\int^{\mi(I)}_{\mi(I)-t}\psi^\ast(u)du.
\label{eq2.1}
\end{eqnarray}
Since $\psi^\ast$ is decreasing we get the inequalities in
(\ref{eq2.1}), while the equality is obvious since
\[
\int^{\mi(I)}_0\psi^\ast(u)du=\int_I\f d\mi.
\]
From (\ref{eq2.1}) it is easily seen that there exists $r\ge0$ such
that $t+r\le\mi(I)$ with
\begin{eqnarray}
\frac{1}{t}\int^{t+r}_r\psi^\ast(u)du=s.  \label{eq2.2}
\end{eqnarray}
It is also easily seen that there exists $E_t$ measurable subset of
$I$ such that
\begin{eqnarray}
\mi(E_t)=t \ \ \text{and} \ \ \int_{E_t}\f
d\mi=\int^{t+r}_r\psi^\ast(u)du  \label{eq2.3}
\end{eqnarray}
since $(X,\mi)$ is non-atomic.

From (\ref{eq2.2}) and (\ref{eq2.3}) we get the conclusion of the
lemma. \qs
\end{Proof}

We now call two measurable subsets of $X$ almost disjoint if
$\mi(A\cap B)=0$.

We give now the following
\begin{Def}\label{def2.1}
A set $\ct$ of measurable subsets of $X$ will be called a tree if
the following conditions are satisfied.
\begin{enumerate}
\item[(i)] $X\in\ct$ and for every $I\in\ct$ we have that
$\mi(I)>0$.
\item[(ii)] For every $I\in\ct$ there corresponds a finite or
countable subset $C(I)\subseteq\ct$ containing at least two
elements such that:
\begin{itemize}
\item[(a)] the elements of $C(I)$ are pairwise almost disjoint
subsets of $I$.
\item[(b)] $I=\cup\, C(I)$.
\end{itemize}
\item[(iii)] $\ct=\bigcup\limits_{m\ge0}\ct_{(m)}$ where
$\ct_0=\{X\}$ and
\[
\ct_{(m+1)}=\bigcup_{I\in\ct_{(m)}}C(I).
\]
\item[(iv)] $\dis\lim_{m\ra+\infty}\sup_{I\in\ct_{(m)}}\mi(I)=0$.
\qs
\end{enumerate}
\end{Def}

From \cite{5} we get the following
\begin{lem}\label{lem2.3}
For every $I\in\ct$ and every $\al$ such that $0<\al<1$ there exists
subfamily $\cf(I)\subseteq Y$ consisting of pairwise almost disjoint
subsets of $I$ such that
\[
\mi\bigg(\bigcup_{J\in \cf(I)}J\bigg)=\sum_{J\in
\cf(I)}\mi(J)=(1-\al)\mi(I).  \text{\qs}
\]
\end{lem}

Let now $(X,\mi)$ be a non-atomic probability measure space and
$\ct$ a tree as in Definition 1.1. We define the associated maximal
operator to the tree $\ct$ as follows: For every  $\f\in L^1(X,\mi)$
and $x\in X$, then
\[
\cm_\ct\f(x)=\sup\bigg\{\frac{1}{\mi(I)}\int_I|\f|d\mi:\;x\in
I\in\ct\bigg\}.
\]
\section{The domain of the extremal problem}  
\noindent

Our aim is to find the exact allowable values of $(f,A,F)$ for which
there exists $\f:(X,\mi)\ra\R^+$ measurable such that
\begin{eqnarray}
\int_X\f d\mi=f,\;\int_X\f^qd\mi=A \ \ \text{and} \ \ |||\f|||_{p,\infty}=F.
\label{eq3.1}
\end{eqnarray}
We find it in the case where $F=1$.

For the beginning assume that $(f,A)$ are such that there exist $\f$
as in (\ref{eq3.1}). We set $g=\f^\ast:[0,1]\ra\R^+$. Then
\[
\int^1_0g=f,\;\int^1_0g^q=A \ \ \text{and} \ \
\big|||g|||_{p,\infty}^{[0,1]}=1
\]
where
\begin{align*}
|||g|||^{[0,1]}_{p,\infty}=\sup\bigg\{&|E|^{-1+\frac{1}{p}}\int_Eg:
\;E\subset[0,1], \ \ \text{Lebesque}\\
&\text{measurable such that } \; |E|>0\bigg\}.
\end{align*}
This is true because of the definition of the decreasing
rearrangement of $\f$ and Lemma \ref{lem2.1}. In fact since $g$ is
decreasing $|||g|||_{p,\infty}$ is equal to
\[
\sup\bigg\{t^{-1+\frac{1}{p}}\int^t_0g: \;\;0<t\le1\bigg\}.
\]
Of course, we should have that $0<f\le1$ and $f^q\le A$. We give now
the following
\begin{Def}\label{def3.1}
If $n\in\N$, and $h:[0,1)\ra\R^+$, $h$ will be called
$\frac{1}{2^n}$-step if it is constant on each interval
\[
\bigg[\frac{i-1}{2^n},\frac{i}{2^n}\bigg), \ \ i=1,2,\ld,2^n.
\text{\qs}
\]
\end{Def}

Now for $n\in\N$ and $0<f\le1$ fixed we set
\begin{align*}
\De_n(f)=\bigg\{&h:\;[0,1]\ra\R^+:\;\;g\;\;\text{is a
$\frac{1}{2^n}$-step function}, \\
&\int^1_0g=f, \ \ |||g|||^{[0,1]}_{p,\infty}\le1\bigg\}.
\end{align*}
Then
\[
\De_n=\De_n(f)\subset L^{p,\infty}([0,1])
\]
where we use the $|||\cdot|||^{[0,1]}_{p,\infty}$ norm for
functions defined on $[0,1]$. $\De_n$ is also convex, that is
\[
h_1,h_2\in\De_n\Rightarrow\frac{h_1+h_2}{2}\in\De_n.
\]
Additionally we have the following
\begin{lem}\label{lem3.1}
$\De_n$ is compact subset of $L^{p,\infty}([0,1])=Y$ where the
topology on $Y$ is that endowed by
$|||\cdot|||^{[0,1]}_{p,\infty}$.
\end{lem}
\begin{Proof}
$(Y,|||\cdot|||_{p,\infty})$ is a Banach space. So, especially a
metric space. So, we just need to prove that $\De_n$ is
sequentially compact.

Let now $(h_i)_i\subset\De_n$. It is now easy to see by a finite
diagonal argument that there exists $(h_{i_j})_j$ subsequence and
$h:[0,1]\ra\R^+$. such that $h_{i_j}\ra h$ uniformly on $[0,1]$.
Then obviously $\int\limits^1_0h=f$,
$|||h|||_{p,\infty}^{[0,1]}\le1$, so $h\in\De_n$. Additionally
\begin{align*}
|||h_{i_j}-h|||^{[0,1]}_{p,\infty}&=\sup\bigg\{|E|^{1+\frac{1}{p}}\int_E|h_{i_j}-h
:\;|E|>0\bigg\}\\
&\le\sup|(h_{i_j}-h)(t)|\;t\in[0,1]
\end{align*}
as $j\ra\infty$. That is
$h_{i_j}\overset{Y}{\longrightarrow}h\in\De_n$. Consequently,
$\De_n$ is a compact subset of $L^{p,\infty}([0,1])$. \qs
\end{Proof}

We give now the following known
\begin{Def}\label{def3.2}
For a closed convex subset $K$ of a topological vector space $Y$,
and for a $y\in K$ we say that $y$ is an extreme point of $K$, if
whenever $y=\frac{x+z}{2}$, with $x,z\in K$ it is implied that
$y=x=z$. We write $y\in ext(K)$.  \qs
\end{Def}
\begin{Def}\label{def3.3}
For a subset $A$ of a topological vector space $Y$ we set
\[
conv(A)=\bigg\{\sum^n_{i=1}\la_ix_i:\;\la_i\ge0,\;x_i\in A,\;
n\in\N^\ast,\;\sum^n_{i=1}\la_i=1\bigg\}.
\]
We call $conv(A)$ the convex hull of $A$.  \qs
\end{Def}

We state now the following well known
\begin{thm}\label{thm3.1}
(Krein\,-\,Milman) Let $K$ be a convex, compact subset of a locally
convex topological vector space $Y$ then
$K=\overline{conv(ext(K))}^Y$ that is $K$ is the closed convex hull
of it's extreme points. \qs
\end{thm}

According now to Lemma \ref{lem3.1} we have that
\[
\De_n=\overline{conv[ext(\De_n)]}^{L^{p,\infty}([0,1])}.
\]
We find now the set $ext(\De_n)$.
\begin{lem}\label{lem3.2}
Let $g\in ext(\De_n)$. Then for every $i\in\{1,2,\ld,2^n\}$ such
that $\Big(\frac{i}{2^n}\Big)^{1-\frac{1}{p}}\le f$ we have that
\[
\sup\bigg\{|E|^{-1+\frac{1}{p}}\int_Eg:|E|=\frac{i}{2^n}\bigg\}=1.
\]
\end{lem}
\begin{Proof}
We prove it first when $i=1$ and
$\Big(\frac{1}{2^n}\Big)^{1-\frac{1}{p}}\le f$. It is now easy to
see that $g\in ext(\De_n)\Leftrightarrow g^\ast\in ext(\De_n)$. So
we just need to prove that
$\int\limits^{1/2^n}_0g^\ast=\Big(\frac{1}{2^n}\Big)^{1-\frac{1}{p}}$.
We write
\[
g^\ast=\sum^{2^n}_{i=1}\al_i\xi_{I_i} \ \ \text{with} \ \
I_i\bigg[\frac{i-1}{2^n},\frac{i}{2^n}\bigg)
\]
and $\al_1\ge\al_{i+1}$ for every $i\in\{1,2,\ld,2^n-1\}$.

Suppose now that $\al_1<2^{n/p}$, and that $\al_1>\al_2$ (the case
$\al_1=\al_2$ is handled in an analogous way).

For a suitable $\e>0$ we set
\[
g_1=\sum^{2^n}_{i=1}\al^{(1)}_i\xi_{I_i}, \ \
g_2=\sum^{2^n}_{i=1}\al^{(2)}_i\xi_{I_i} \ \ \text{where} \ \
\left.\begin{array}{cc}
  \al^{(1)}_1=\al_1+\e, & \al^{(1)}_2=\al_2-\e \\
  \al^{(2)}_1=\al_1-\e, & \al^{(2)}_2=\al_2+\e \\
\end{array}\right\}
\]
and $\al^{(1)}_k=\al^{(2)}_k=\al_k$ for every $k>2$.

Since $\al_1<2^{n/p}$ we can find small enough $\e>0$ such that
$g_i$ satisfy $|||g_i|||^{[0,1]}_{p,\infty}\le1$, for $i=1,2$.
Indeed, for $i=1$, we need to prove that for small enough $\e>0$
\begin{eqnarray}
\int^t_0g_1\le t^{1-\frac{1}{p} }  \label{eq3.2}
\end{eqnarray}
for every $t\in[0,1)$, since $g_1$ is decreasing.

(\ref{eq3.2}) is now obviously true for $t\ge\frac{2}{2^n}$ since
\begin{eqnarray}
\int^t_0g_1=\int^t_0g^\ast \ \ \text{for every such} \ \ t
\label{eq3.3}
\end{eqnarray}
(\ref{eq3.2}) is also true for $t=0,\frac{1}{2^n}$. But then it
remains true for every $t\in\Big(0,\frac{1}{2^n}\Big)$ since the
function $t\mapsto\int\limits^t_0g_1$ represents a straight line on
$\Big[0,\frac{1}{2^n}\Big]$ and $t^{1-\frac{1}{}p}$ is concave
there. Analogously for the interval
$\Big[\frac{1}{2^n},\frac{2}{2^n}\Big]$. That is we proved
$|||g_1|||^{[0,1]}_{p,\infty}\le1$.

Obviously, $\int\limits^1_0g_i=f$, so that $g_i\in\De_n$, for
$i=1,2$. But $g^\ast=\frac{g_1+g_2}{2}$, with $g_i\neq g$ and
$g_i\in\De_n$, $i=1,2,$, a contradiction since $g^\ast\in
ext(\De_n)$. So,
\[
\al_1=2^{n/p} \ \ \text{and} \ \ \int^{1/2}_0g^\ast=
\bigg(\frac{1}{2^n}\bigg)^{1-\frac{1}{p}},
\]
what we wanted to prove. In the same way we prove that for
$i\in\{1,2,\ld,2^n\}$ such that
\[
\bigg(\frac{i+1}{2^n}\bigg)^{1-\frac{1}{p}}\le f, \ \ \text{if} \
\ \int^{i/2^n}_0g^\ast=\bigg(\frac{i}{2^n}\bigg)^{1-\frac{1}{p}} \; \text{then}
\; \int^{i+1/2^n}_0g^\ast=\bigg(\frac{i+1}{2^n}\bigg)^{1-\frac{1}{p}}.
\]
The lemma is now proved.  \qs
\end{Proof}

Let now $g\in ext(\De_n)$ and
$k=\max\Big\{i\le2^n:\Big(\frac{i}{2^n}\Big)^{1-\frac{1}{p}}\le
f\Big\}$, so if we suppose that $f<1$ we have that
\[
\bigg(\frac{k}{2^n}\bigg)^{1-\frac{1}{p}}\le
f<\bigg(\frac{k+1}{2^n}\bigg)^{1-\frac{1}{p}}.
\]
By Lemma \ref{lem3.2}
\[
\int^{k/2^n}_0g^\ast=\bigg(\frac{k}{2^n}\bigg)^{1-\frac{1}{p}}.
\]
But by using the reasoning of the previous lemma it is easy to see
that
\[
\int^{k+1/2^n}_0g^\ast=f,
\]
which gives
\[
\int^{k+1/2^n}_{k/2^n}g^\ast=f-\bigg(\frac{k}{2^n}\bigg)^{1-\frac{1}{p}}
\Rightarrow\al_{k+1}=2^n\cdot f-2^{n/p}\cdot k^{1-\frac{1}{p}}.
\]
Additionally $\al_i=0$ for $i>k+1$.

From the above we obtain the following
\begin{cor}\label{cor3.1}
Let $g\in ext(\De_n)$. Then
$g^\ast=\sum\limits^{2^n}_{i=1}\al_i\xi_{I_i}$ where
\[
\al_i=2^{n/p}\Big(i^{1-\frac{1}{p}}-(i-1)^{1-\frac{1}{p}}\Big) \ \
\text{for} \ \ i=1,2,\ld,k
\]
and
\[
\al_{k+1}=2^nf-2^{n/p}\cdot k^{1-\frac{1}{p}}, \ \ \al_i=0, \ \ i>k+1,
\]
where
\[
k=max\bigg\{i\le2^n:\bigg(\frac{i}{2^n}\bigg)^{1-\frac{1}{p}}\le
f\bigg\}.  \text{\qs}
\]
\end{cor}
\begin{rem}\label{rem3.1}
Actually it is easy to see that the above functions described in
Corollary \ref{cor3.1} are exactly the extreme points of $\De_n$.
\qs
\end{rem}

We estimate now the $L^q$-norm of every $g\in ext(\De_n)$.

We state it as
\begin{lem}\label{lem3.3}
Let $g\in ext(\De_n)$ and $A=\int\limits^1_0g^q$, then $A\le\Ga
f^{p-q/p-1}+\ce_n(f)$ where
\[
\Ga=\bigg(\frac{p-1}{p}\bigg)^q\frac{p}{p-q} \ \ \text{and} \ \
\ce_n(f)=\frac{\al^q_{k+1}}{2^n}=\frac{(2^nf-2^{n/p}k^{1-\frac{1}{p}})^q}
{2^n}.
\]
\end{lem}
\begin{Proof}
For $g$ we write $g^\ast=\ssum^{2^n}_{i=1}\al_i\xi_{I_i}$, where
$\al_i$ are given in Corollary \ref{cor3.1}. Then
\begin{eqnarray}
A=\int^1_0(g^\ast)^q=\bigg[\bigg(\sum^k_{i=1}\al^q_i\bigg)+\al^q_{k+1}\bigg]\cdot
\frac{1}{2^n}. \label{eq3.4}
\end{eqnarray}
Now for $i\in\{1,2,\ld,k\}$
\begin{align}
\al^q_i&=\bigg[2^{n/p}\Big(i^{1-\frac{1}{p}}-(i-1)^{1-\frac{1}{p}}\Big)\bigg]^q
=\bigg\{2^n\bigg[\bigg(\frac{i}{2^n}\bigg)^{1-\frac{1}{p}}-\bigg(
\frac{i-1}{2^n}\bigg)^{1-\frac{1}{p}}\bigg]\bigg\}^q \nonumber \\
&=\bigg[2^n\int^{i/2^n}_{i-1/2^n}\psi\bigg]^q  \label{eq3.5}
\end{align}
where $\psi:(0,1]\ra\R^+$ is defined by
$\psi(t)=\frac{p-1}{p}t^{-1/p}$. By (\ref{eq3.5}) and in view of
Holder's inequality we have that for $i\in\{1,2,\ld,k\}$
\begin{eqnarray}
\al^q_i\le2^n\int^{i/2^n}_{i-1/2^n}\psi^q.  \label{eq3.6}
\end{eqnarray}
Summing up relations (\ref{eq3.6}) we have that
\begin{eqnarray}
\sum^k_{i=1}\al^q_i\le2^n\int^{k/2^n}_0\psi^q=2^n\cdot\Ga\cdot\bigg(\frac{k}{2^n}
\bigg)^{1-\frac{q}{p}}.  \label{eq3.7}
\end{eqnarray}
Additionally from the definition of $k$ we have that
\begin{eqnarray}
\bigg(\frac{k}{2^n}\bigg)^{1-\frac{1}{p}}\le f \;\Rightarrow\;
k^{1-\frac{q}{p}}\le(2^n)^{1-\frac{q}{p}}\cdot f^{p-q/p-1}.
\label{eq3.8}
\end{eqnarray}
From (\ref{eq3.4}), (\ref{eq3.7}) and (\ref{eq3.8}) we obtain
\[
A\le\bigg[2^n\cdot\Ga\cdot
f^{p-q/p-1}+\al^q_{k+1}\bigg]\frac{1}{2^n} =\Ga
f^{p-q/p-1}+\ce_n(f)
\]
and Lemma \ref{lem3.3} is proved. \qs
\end{Proof}
\begin{cor}\label{cor3.2}
For every $g\in\De_n$
\[
A\le\Ga f^{p-q/p-1}+\ce_n(f), \ \ \text{where} \ \
A=\int^1_0g^q.
\]
\end{cor}
\begin{Proof}
This is true, of course, for $g\in ext(\De_n)$, and so also for
$g\in conv(ext\De_n)$, since $t\mapsto t^q$ is convex for $q>1$ on
$\R^+$. It remains true for
$g\in\overline{conv(ext(\De_n))}^{L^{p,\infty}([0,1])}$ using a
simple continuity argument. In fact we just need the continuity of
the identity operator if it is viewed as: $I:L^{p,\infty}([0,1])\ra
L^q([0,1])$. See \cite{4}.

Using now Krein\,-\,Milman Theorem the Corollary is proved. \qs
\end{Proof}

We have now the following
\begin{cor}\label{cor3.3}
Let $\f:(X,\mi)\ra\R^+$ such that
\[
\int_X\f d\mi=f, \ \ \int_X\f^qd\mi=A, \ \
|||\f|||_{p,\infty}\le1.
\]
Then
\[
f^q\le A\le\Ga f^{p-q/p-1}.
\]
\end{cor}
\begin{Proof}
Let $g=\f^\ast:[0,1]\ra\R^+$. There exist $\f_n\frac{1}{2^n}$-simple
functions, for every $n$ such that $g_n\le g_{n+1}\le g$ and $g_n$
converges almost everywhere to $g$. But then by defining
\[
f_n=\int^1_0\f_n, \ \ A_n=\int^1_0\f^q_n
\]
we have that
\begin{eqnarray}
g_n\in\De_n(f_n) \ \ \text{so that} \ \ A_n\le\Ga f^{p-q/p-1}_n+
\ce_n(f_n).  \label{eq3.9}
\end{eqnarray}
By the monotone convergence theorem $f_n\ra f$, $A_n\ra A$. Moreover
\[
\ce_n(f_n)=\frac{(2^nf_n-k_n^{1-\frac{1}{p}}2^{n/p})^q}{2^n}
\]
where $k_n$ satisfy
\[
\bigg(\frac{k_n}{2^n}\bigg)^{1-\frac{q}{p}}\le
f_n<\bigg(\frac{k_n+1}{2^n}\bigg)^{1-\frac{1}{p}}.
\]
As a consequence
\begin{align*}
\ce_n(f_n)&=(2^n)^{q-1}\bigg[f_n-\bigg(\frac{k_n}{2^n}\bigg)^{1-\frac{1}{p}}\bigg]^q<
(2^n)^{q-1}\bigg[\bigg(\frac{k_n+1}{2^n}\bigg)^{1-\frac{1}{p}}-\bigg(\frac{k_n}{2^n}\bigg)
^{1-\frac{1}{p}}\bigg]^q \\
&\le(2^n)^{q-1}\bigg[\bigg(\frac{1}{2^n}\bigg)^{1-\frac{1}{q}}\bigg]^q=
\bigg(\frac{1}{2^{1-\frac{q}{p}}}\bigg)^n\ra0, \ \ \text{as} \ \
n\ra\infty
\end{align*}
where in the second inequality we used the known
\[
(t+s)^\al\le t^\al+s^\al \ \ \text{for} \ \ t,s\ge0, \ \ 0<\al<1.
\]
Now (\ref{eq3.9}) gives the corollary.  \qs
\end{Proof}

In fact the converse of Corollary \ref{cor3.3} is also true.
\begin{thm}\label{thm3.2}
For $0<f\le1$, $A>0$ the following are equivalent
\vspace*{0.2cm} \\
i)\;$f^q\le A\le\Ga f^{p-q/p-1}$ \vspace*{0.2cm} \\
ii)\; $\exists\;\f:(X,\mi)\ra\R^+$ such that
\[
\int_X\f d\mi=f, \ \ \int_X\f^qd\mi=A, \ \
|||\f|||_{p,\infty}\le1. \text{\qs}
\]
\end{thm}

We prove first the following
\begin{lem}\label{lem3.4}
Let $\al\in(0,1)$ and $(f,A)$ such that
\begin{eqnarray}
f\lneq\al^{1-\frac{1}{p}} \label{eq3.10}
\end{eqnarray}
\begin{eqnarray}
f^q\lneq\al^{q-1}A \label{eq3.11}
\end{eqnarray}
\begin{eqnarray}
A\le\Ga f^{p-q/p-1}A. \label{eq3.12}
\end{eqnarray}

Then there exists $g:[0,\al]\ra\R^+$ such that
\[
\int^\al_0g=f, \ \ \int^\al_0g^q=A, \ \ \text{and} \ \
|||g|||^{[0,\al]}_{p,\infty}=1
\]
where
\[
|||g|||^{[0,\al]}_{p,\infty}=\sup\left\{\begin{array}{ll}
                                           & E\ \ \text{measurable subset
of} \ \ [0,\al] \\ [-2ex]
                     |E|^{-1+\frac{1}{p}}\int_Eg:                      &  \\ [-2ex]
                                           & \text{such that} \; |E|>0.
                                        \end{array}\right\}
\]
\end{lem}
\begin{Proof}
We search for a $g$ of the form
\[
g:=\left\{\begin{array}{cc}
  \frac{p-1}{p}t^{-1/p}, & 0<t\le c_1 \\
  \mi_2, & c_1<t\le\al \\
\end{array}\right.
\]
for suitable constant $c_1\mi_2$.

We must have that
\begin{eqnarray}
\int^\al_0g=f\;\Leftrightarrow\;c_1^{1-\frac{1}{p}}+\mi_2(\al-c_1)=f.
\label{eq3.13}
\end{eqnarray}
Additionally $g$ must satisfy
\begin{eqnarray}
\int^\al_0g^q=A\;\Leftrightarrow\;\Ga
c_1^{1-\frac{q}{p}}+\mi^q_2(\al-c_1)=A.  \label{eq3.14}
\end{eqnarray}
(\ref{eq3.13}) gives
\begin{eqnarray}
\mi_2=\frac{f-c^{1-\frac{1}{p}}}{a-c_1}  \label{3.15}
\end{eqnarray}
so (\ref{eq3.14}) becomes
\begin{eqnarray}
\Ga
c^{1-\frac{q}{p}}_1+\frac{(f-c_1^{1-\frac{1}{P}})^q}{(\al-c_1)^{q-1}}=A.
\label{eq3.16}
\end{eqnarray}
We we search for a $c_1\in(0,\al)$ such that
\[
T(c_1)=A \ \ \text{where} \ \ T:[0,\al)\ra\R^+
\]
defined by
\[
T(t)=\Ga
t^{1-\frac{q}{p}}+\frac{(f-t^{1-\frac{1}{p}})^q}{(\al-t)^{q-1}}.
\]

Observe that $T(0)=\frac{f^q}{\al^{q-1}}\lneq A$ because of
(\ref{eq3.11}) and that $T(f^{p/p-1})=\Ga f^{p-q/p-1}\ge A$. Now
because of the continuity of $T$, we have that there exists
$c_1\in(0,f^{p/p-1}]$ such that $T(c_1)=A$. Then $c_1\in(0,\al)$
because of (\ref{eq3.10}), and if we define $\mi_2$ by
(\ref{eq3.15}), we guarantee (\ref{eq3.13}) and (\ref{eq3.14}).

We need to prove now that $|||g|||^{[0,\al]}_{p,\infty}=1$.

Obviously, because of the form of $g$,
$|||g|||^{[0,\al]}_{p,\infty}\ge1$. So we have to prove that
\begin{eqnarray}
\int^t_0g\le t^{1-\frac{1}{p}}, \ \ \fa\;t\in(0,\al].
\label{eq3.17}
\end{eqnarray}
This is of course true for $t\in[0,c_1]$. For $t\in(c_1,\al]$
\[
\int^t_0g=c^{1-\frac{1}{p}}_1+\mi_2(t-c_1)=:G(t).
\]
Since $G(c_1)=c_1^{1-\frac{1}{p}}$, $G(\al)=f<\al^{1-\frac{1}{p}}$
and $t\mapsto t^{1-\frac{1}{p}}$ is concave on $(c_1,\al]$
(\ref{eq3.17}) is true. Thus Lemma \ref{lem3.4} is proved.  \qs
\end{Proof}

We have now the\vspace*{0.2cm} \\
{\bf Proof of Theorem 3.2:} We have to prove the direction
i)\;$\Rightarrow$\;ii).

Indeed if $f^q\lneq A\le\Ga f^{p-q/p-1}$ and $f<1$ we apply Lemma
\ref{lem3.4}.

If $f^q=A$, with $0<f\le1$ we set $g$ by $g(t)=f$, for every
$t\in[0,1]$ while if $f=1\le A\le\Ga$ a simple modification of Lemma
\ref{lem3.4} gives the result. \qs \medskip

We conclude Section 3 with the following theorem which can be
proved easily using all the above.
\begin{thm}\label{thm3.3}
For $f,A$ such that $f<1$, $A>0$ the following are
equivalent:\vspace*{0.2cm} \\
i)\;$f^q\lneq A\le\Ga f^{p-q/p-1}$ \vspace*{0.2cm} \\
ii)\;$\exists\;\f:(X,\mi)\ra\R^+$ such that
\[
\int_X\f d\mi=f, \ \ \int_X\f^qd\mi=A, \ \
|||\f|||_{p,\infty}=1.  \text{\qs}
\]
\end{thm}
\begin{rem}\label{rem3.2}
Theorem \ref{thm3.3} is completed if we mention that for $f=1$ the
following are equivalent: \vspace*{0.2cm} \\
i)\;$f=1\le A\le\Ga$ \vspace*{0.2cm} \\
ii)\;$\exists\;\f:(X,\mi)\ra\R^+$ such that $\int\limits_X\f
d\mi=1$, $\int\limits_X\f^qd\mi=A$, $|||\f|||_{p,\infty}=1$. \qs
\end{rem}
\section{The Extremal Problem} 
\noindent

Let $\cm_\ct=\cm$ the dyadic maximal operator associated to the tree
$\ct$, on the probability non-atomic measure space $(X,\mi)$.

Our aim is to find
\begin{align*}
T_{f,A,F}(\la)=\sup\bigg\{\mi(\{\cm\f\ge\la\}):\f\ge0,\;
\int_X\f d\mi=f,\;&\int_Z\f^qd\mi=A, \\
&|||\f|||_{p,\infty}=F\bigg\}
\end{align*}
for all the allowable values of $f,A,F$.

We find it in the case where $F=1$.

We write $T_{f,A}(\la)$ for $T_{f,A,1}(\la)$.

In order to find $T_{f,A}(\la)$ we find first the following
\begin{align}
T_{f,A}^{(1)}(\la)=\sup\bigg\{\mi(\{\cm\f\ge\la\}):\f\ge0,\;
\int_X\f d\mi=f,;&\int_X\f^qd\mi=A,\nonumber \\
&|||\f|||_{p,\infty}\le1\bigg\}. \label{eq4.1}
\end{align}
The domain of this extremal problem is the following:
\[
D=\Big\{(f,A):\;0<f\le1, \ \ f^q\le A\le\Ga f^{p-q/p-1}\Big\}.
\]
Obviously, $T^{(1)}_{f,A}(\la)=1$, for $\la\le f$.

Now for $\la>f$ and $(f,A)\in D$.

Let $\f$ as in (\ref{eq4.1}). Consider the decreasing rearrangement
of $\f$, $g=\f^\ast:[0,1]\ra\R^+$. Then
\[
\int^1_0g=f,\ \ \int^1_0g^q=A, \ \
|||g|||^{[0,1]}_{p,\infty}\le1.
\]
Consider also $E=\{\cm\f\ge\la\}\subseteq X$.

Then $E$ is the almost disjoint union of elements of $\ct$, let
$(I_j)_j$. In fact we just need to consider the elements $I$ of
$\ct$, maximal under the condition
\begin{eqnarray}
\frac{1}{\mi(I)}\int_I\f d\mi\ge\la.  \label{eq4.2}
\end{eqnarray}
We, then, have $E=\bigcup\limits_jI_j$ and $\int\limits_E\f
d\mi\ge\la\mi(E)$ because of (\ref{eq4.2}). Then according to Lemma
\ref{lem2.1} we have that $\int\limits^\al_0g\ge\al\la$ where
$\al=\mi(E)$. That is we proved that
\begin{eqnarray}
T^{(1)}_{f,A}(\la)\le\De_{f,A}(\la)  \label{eq4.3}
\end{eqnarray}
where
\begin{align}
\De_{f,A}(\la)=\sup\bigg\{&\al\in(0,1]:\;\exists\;g:[0,1]\ra\R^+
\ \ \text{with}
\int^1_0g=f,\;\int^1_0g^q=A,\nonumber\\
&|||g|||^{[0,1]}_{p,\infty}\le1 \ \ \text{and} \ \
\int^\al_0g\ge\al\la\bigg\}.  \label{eq4.4}
\end{align}
We prove now the converse inequality in (\ref{eq4.3}) by proving the
following
\begin{lem}\label{lem4.1}
Let $g$ be as in (\ref{eq4.4}) for a fixed $\al\in(0,1]$. Then there
exists $\f:(X,\mi)\ra\R^+$ such that
\[
\int_X\f d\mi=f, \ \ \int_X\f^qd\mi=A, \ \
|||\f|||_{p,\infty}\le1 \ \ \text{and} \ \
\mi(\{\cm\f\ge\la\})\ge\al.
\]
\end{lem}
\begin{Proof}
Lemma \ref{lem2.3} guarantees the existence of a sequence
$(I_j)_j$ of pairwise almost disjoint elements of $\ct$ such that
\begin{eqnarray}
\mi(\cup I_j)=\sum\mi(I_j)=\al.  \label{eq4.5}
\end{eqnarray}
Consider now the finite measure space $([0,\al],|\cdot|)$ where
$|\cdot|$ is the Lebesque measure. Then since $\int\limits^\al_0
g\ge\al\la$ and (\ref{eq4.5}) holds, applying Lemma \ref{lem2.2}
repeatedly, we obtain the existence of a sequence $(A_j)$ of
Lebesque measurable subsets of $[0,\al]$ such that the following
hold:
\[
(A_j)_j \; \text{is a pairwise disjoint family}, \;
\cup A_j=[0,\al], \; |A_j|=\mi(I_j), \;
\frac{1}{|A_j|}\int_{A_j}g\ge\la.
\]
Then we define $g_j:[0,|A_j|]\ra\R^+$ by $g_j=(g/A_j)^\ast$. Define
also for every $j$ a measurable function $\f_j:I_j\ra\R^+$ so that
$\f^\ast_j=g_j$. The existence of such a function is guaranteed by
the fact that $(I_j,\mi/I_j)$ is non-atomic. Here we mean
\[
\mi/I_j(A)=\mi(A\cap I_j) \ \ \text{for every} \ \
A\subseteq I_j.
\]
Since $(I_j)$ is almost pairwise disjoint family we produce a
$\f^{(1)}:\cup I_j\ra\R^+$ measurable such that $\f^{(1)}/I_j=\f_j$.
We set now $Y=X\sm\cup I_j$ and $h:[0,1-\al]\ra\R^+$ by
$h=(g/[\al,1])^\ast$. Then since $\mi(Y)=1-\al$ there exists
$\f^{(2)}:Y\ra\R^+$ such that $(\f^{(2)})^\ast=h$.

Set now $\f=\left\{\begin{array}{ccc}
  \f^{(1)}, & \text{on} & \cup I_j \\
  \f^{(2)}, & \text{on} & Y. \\
\end{array}\right.$

It is easy to see from the above construction that $\int\limits_X\f
d\mi=f$, $\int\limits_X\f^qd\mi=A$ and $|||\f|||_{p,\infty}\le1$.

Additionally
\[
\frac{1}{\mi|I_j|}\int_{I_j}\f
d\mi=\frac{1}{|A_j|}\int_{A_j}g\ge\la \ \ \text{for every} \ \
j
\]
that is
\[
\{\cm\f\ge\la\}\supseteq\cup I_j, \ \ \text{so} \ \
\mi(\{\cm\f\ge\la\})\ge\al
\]
and the lemma is proved.  \qs
\end{Proof}

It is now not difficult to see that we can replace the inequality
$\int\limits^\al_0g\ge\al\la$ in the definition of $\De_{f,A}(\la)$
by equality, thus giving $S_{f,A}(\la)$, in such a way that
(\ref{eq4.3}) remains true, that is
\begin{eqnarray}
T^{(1)}_{f,A}(\la)=\De_{f,A}(\la)=S_{f,A}(\la).
\label{eq4.6}
\end{eqnarray}
This is true since if $g$ is as in (\ref{eq4.4}) there exists
$\bi\ge\al$ such that $\int\limits^\bi_0g=\bi\la$.

For $(f,A)\in D$ we set
\[
G_{f,A}(\la)=\sup\bigg\{\mi(\{\cm\f\ge\la\}):\;\f\ge0,\;\int\limits_X\f
d\mi=f,\; \int\limits_X\f^qd\mi=A\bigg\}.
\]
It is obvious that $T^{(1)}_{f,A}(\la)\le G_{f,A}(\la)$.

As a matter of fact $G_{f,A}(\la)$ has been computed in \cite{3} and
was found to be
\begin{eqnarray}
G_{f,A}(\la)=\left\{\begin{array}{cl}
  1, & \la\le f \\
  \frac{f}{\la}, & f<\la<\Big(\frac{A}{f}\Big)^{1/q-1} \\
  k, & \Big(\frac{A}{f}\Big)^{1/q-1}\le\la \\
\end{array}\right.  \label{eq4.7}
\end{eqnarray}
where $k$ is the unique root of the equation
\[
\frac{(f-\al\la)^q}{(1-\al)^{q-1}}+\al\la^q=A \ \ \text{on} \
\ \bigg[0,\frac{f}{\la}\bigg], \ \ \text{when} \ \
\la>\bigg(\frac{A}{f}\bigg)^{1/q-1}.
\]
We have now the following
\begin{prop}\label{prop4.1}
For $(f,A)\in D$, then
\[
T^{(1)}_{f,A}(\la)\le\min\bigg\{1,G_{f,A}(\la),\frac{1}{\la^p}\bigg\}.
\]
\end{prop}
\begin{Proof}
We just need to see that $\mi(\{\cm\f\ge\la\})\le\frac{1}{\la^p}$
for every $\f$ such that $|||\f|||_{p,\infty}\le1$. But if
$E=\{\cm\f\ge\la\}$ we have by the definition of the norm
$|||\cdot|||_{p,\infty}$ that
$\int\limits_E\cm\f\le\mi(E)^{1-\frac{1}{p}}$. But by (\ref{eq1.3})
$\int\limits_E\cm\f\ge\la\mi(E)$, so that
\[
\la\mi(E)\le\mi(E)^{1-\frac{1}{p}}\Rightarrow\mi(E)\le\frac{1}{\la^p}.
\]
So Proposition \ref{prop4.1} is true.  \qs
\end{Proof}

We prove now the converse of Proposition \ref{prop4.1} in three
steps.
\begin{prop}\label{prop4.2}
Let $(f,A)\in D$ and $\la$ such that
\begin{eqnarray}
\frac{f}{\la}=\min\bigg\{1,G_{f,A}(\la),\frac{1}{\la^p}\bigg\}.
\label{eq4.8}
\end{eqnarray}
Then $T^{(1)}_{f,A}(\la)=\frac{f}{\la}$.
\end{prop}
\begin{Proof}
We use Lemma \ref{lem3.4} and equations (\ref{eq4.6}). Because of
(\ref{eq4.6}) we need to find $g:[0,1]\ra\R^+$ such that
\[
\int^1_0g=f, \ \ \int^1_0g^q=A, \ \ |||g|||_{p,\infty}\le1 \ \
\text{and} \ \ \int^{f/\la}_0g=\frac{f}{\la}\cdot\la=f
\]
that is $g$ should be defined on $[0,f/\la]$.

We apply Lemma \ref{lem3.4}, with $\al=\frac{f}{\la}$.

In fact, since (\ref{eq4.8}), is true we have that
$G_{f,A}(\la)=\frac{f}{\la}$ so, $\la<\Big(\frac{A}{f}\Big)^{1/q-1}$
which gives (\ref{eq3.11}), while $\frac{f}{\la}\le\frac{1}{\la^p}$
gives (\ref{eq3.10}). In fact Lemma \ref{lem3.4} works even with
equality on (\ref{eq3.10}) as it is easily can be seen. So, in view
of (\ref{eq4.6}) we have $T^{(1)}_{f,A}(\la)\ge f/\la$ and the
proposition is proved. \qs
\end{Proof}

At the next step we have
\begin{prop}\label{prop4.3}
Let $(f,A)\in D$ and $\la$ such that
\begin{eqnarray}
k=\min\bigg\{1,G_{f,A}(\la)\frac{1}{\la^p}\bigg\}. \label{eq4.9}
\end{eqnarray}
Then $T^{(1)}_{f,A}(\la)=k$.
\end{prop}
\begin{Proof}
Obviously (\ref{eq4.9}) gives $\la\ge\Big(\frac{A}{f}\Big)^{1/q-1}$.

We prove that there exists $g:[0,1]\ra\R^+$ such that
\begin{eqnarray}
\int^k_0g=k\la,  \ \ \int^1_0g=f, \ \ \int^1_0g^q=A \ \
\text{and} \ \ |||g|||_{p,\infty}\le1. \label{eq4.10}
\end{eqnarray}
For this purpose we define:
\[
g:=\left\{\begin{array}{ccc}
  \la, & \text{on} & [0,k] \\
  \frac{f-k\la}{1-k}, & \text{on} & (k,1]. \\
\end{array}\right.
\]
Then, obviously, the first two conditions in (\ref{eq4.10}) are
satisfied, while
\[
\int^1_0g^q=\frac{(f-k\la)^q}{(1-k)^{q-1}}+k\la^q=A,
\]
by the definition of $k$.

Moreover $|||g|||_{p,\infty}\le1$. This is true since $k\la\le
k^{1-\frac{q}{p}}$, $f\le1$ and the fact that $g$ is constant on
each of the intervals $[0,k]$ and $(k,1]$. So we proved that
$T^{(1)}_{f,A}(\la)\ge k$, that is what we wanted to prove.  \qs
\end{Proof}

At last we prove
\begin{prop}\label{prop4.4}
Let $(f,A)\in D$ and $\la$ such that
\begin{eqnarray}
\frac{1}{\la^p}=\min\bigg\{1,G_{f,A}(\la),\frac{1}{\la^p}\bigg\}.
\label{eq4.11}
\end{eqnarray}
Then $T^{(1)}_{f,A}(\la)=\frac{1}{\la^p}$.
\end{prop}
\begin{Proof}
As before we search for a function $g$ such that
\begin{equation}
\int^1_0g=f, \ \ \int^1_0g^q=A, \ \ |||g|||_{p,\infty}\le1 \ \
\text{and} \ \
\int^{1/\la^p}_0g=\frac{1}{\la^p}\cdot\la=\frac{1}{\la^{p-1}}.
\label{eq4.12}
\end{equation}
We define
\[
\vthi_\la=\frac{\Ga}{\la^{p-q}}+\frac{\Big(f-\frac{1}{\la^{p-1}}\Big)^q}
{\Big(1-\frac{1}{\la^p}\Big)^{q-1}},
\]
and we consider two cases:
\begin{enumerate}
\item[i)]\;$\vthi_\la>A$
\end{enumerate}

We search for a function of the form
\begin{eqnarray}
g:=\left\{\begin{array}{cc}
  \Big(1-\frac{1}{p}\Big)t^{-1/p}, & 0<t\le c_1 \\
  \mi_2, & c_1<t\le\frac{1}{\la^p} \\
  \mi_3, & \frac{1}{\la^p}<t<1 \\
\end{array}\right.  \label{eq4.13}
\end{eqnarray}
for suitable constants $c_1\le\frac{1}{\la^p}$, $\mi_2,\mi_3$.
Then in view of (\ref{eq4.12}) the following must hold:
\begin{eqnarray}
c^{1-\frac{1}{p}}_1+\mi_2\bigg(\frac{1}{\la^p}-c_1\bigg)=\frac{1}{\la^{p-1}}
\label{eq4.14}
\end{eqnarray}
\begin{eqnarray}
c^{1-\frac{1}{p}}_1+\mi_2\bigg(\frac{1}{\la^p}-c_1\bigg)+\mi_3\bigg(1-\frac{1}
{\la^p}\bigg)=f  \label{eq4.15}
\end{eqnarray}
\begin{eqnarray}
\Ga
c^{1-\frac{q}{p}}_1+\mi_2^q\bigg(\frac{1}{\la^p}-c_1\bigg)+\mi_3^q
\bigg(1-\frac{1}{\la^p}\bigg)=A.  \label{eq4.16}
\end{eqnarray}
Notice that the condition $|||g|||_{p,\infty}\le1$ is
automatically satisfied because of the form of $g$ and the
previous stated relations.

Now (\ref{eq4.14}) and (\ref{eq4.15}) give
\begin{eqnarray}
\mi_3=\frac{f-\frac{1}{\la^{p-1}}}{1-\frac{1}{\la^p}},
\label{eq4.17}
\end{eqnarray}
and
\begin{eqnarray}
\mi_2=\frac{\frac{1}{\la^{p-1}}-c^{1-\frac{1}{p}}}{\frac{1}{\la^p}-c_1},
\label{eq4.18}
\end{eqnarray}
while (\ref{eq4.16}) gives $T(c_1)=A$ where $T$ is defined on
$\Big[0,\frac{1}{\la^p}\Big)$ by
\[
T(c)=\Ga
c^{1-\frac{q}{p}}+\frac{\Big(\frac{1}{\la^{p-1}}-c_1^{1-\frac{1}{p}}\Big)^q}
{\Big(\frac{1}{\la^p}-c\Big)^{q-1}}+\frac{\Big(f-\frac{1}{\la^{p-1}}\Big)^q}
{\Big(1-\frac{1}{\la^p}\Big)^{q-1}}.
\]
Then
\[
T(0)=\frac{1}{\la^{p-q}}+\frac{\Big(f-\frac{1}{\la^{p-1}}\Big)^q}
{\Big(1-\frac{1}{\la^p}\Big)^{q-1}}.
\]
It is now easy to see that $T(0)\le A$ by using that
$F:[0,f/\la]\ra\R^+$ defined by
\[
F(t)=\frac{(f-t\la)^q}{(1-t)^{q-1}}+t\la^q
\]
is increasing, and the definition of $G_{f,A}(\la)$.

Moreover $\dis\lim_{c\ra\frac{1^-}{\la^p}}T(c)=\vthi_\la>A$, so by
continuity of the function $t$, we end case i). Now for
\begin{enumerate}
\item[ii)] \; $\vthi\la\le A$ we search for a function of the
form
\end{enumerate}
\[
g:=\left\{\begin{array}{cc}
  \Big(1-\frac{1}{p}\Big)^{t^{-1/p}}, & 0<t\le c_1 \\
  \mi_2, & c_1<t\le1 \\
\end{array}\right.
\]
where $\frac{1}{\la^p}<c_1$. Similar arguments as in case i) give
the result.  \qs
\end{Proof}

From Propositions \ref{prop4.1}\,-\,\ref{prop4.4} we have now of
course
\begin{thm}\label{thm4.1}
For
\[
(f,A)\in D, \ \
T^{(1)}_{f,A}(\la)=\min\bigg\{1,G_{f,A}(\la),\frac{1}{\la^p}\bigg\}.
\text{\qs}
\]
\end{thm}
\begin{rem}\label{rem4.1}
Notice that $T_{f,A}(\la)=T^{(1)}_{f,A}(\la)$ for every $f,A$ such
that $f^q<A\le\Ga f^{p-q/p-1}$ and $0<f\le1$. Indeed suppose that
$\al=T^{(1)}_{f,A}(\la)$. Then there exists $g:[0,1]\ra\R^+$ such
that
\begin{eqnarray}
\int^1_0g=f, \ \ \int^1_0g^q=A, \ \ \int^\al_0g=\al\la \ \
\text{and} \ \ |||g|||_{p,\infty}\le1.  \label{eq4.19}
\end{eqnarray}
It is easy to see that for every $\e>0$, small enough we can
produce from $g$ a function $g_\e$ satisfying
\[
\int^{\al-\e}_0g_\e\ge(\al-\e)\la, \ \ \int^1_0g_\e=f, \ \
\int^1_0g_\e=A+\de_\e \ \ \text{and} \ \ |||g_\e|||_{p,\infty}=1
\]
and $\dis\lim_{\e\ra0^+}\de\e=0$. This and continuity reasons shows
$T_{f,A}(\la)=\al$.
\begin{enumerate}
\item[iii)] The case $A=f^q$ can be worked out separately
because there is essentially unique function $g$ satisfying
$\int\limits^1_0g=f$, $\int\limits^1_0g^q=f^q$, namely the
constant function with value $f$.  \qs
\end{enumerate}
\end{rem}

Scaling all the above we have that
\begin{thm}\label{thm4.2}
For $f,A$ such that $f^q<A\le\Ga f^{p-q/p-1}F^{p(q-1)/p_1}$ and
$0<f\le F$ the following hold
\begin{align}
\sup\bigg\{\mi(\{\cm\f\ge\la\}):\;\f\ge0,\;\int_X\f d\mi=f,\;
&\int_X\f^qd\mi=A,\;|||\f|||_{p,\infty}=F\bigg\} \nonumber \\
&=\min\bigg\{1,G_{f,A}(\la),\frac{F^p}{\la^p}\bigg\}
\label{eq4.20}
\end{align}
and
\[
\sup\bigg\{\|\cm\f\|_{p,\infty}:\;\f\ge0,\;\int_X\f
d\mi=f,\;\int_X\f^qd\mi=A, \;|||\f|||_{p,\infty}=F\bigg\}=F.
\text{\qs}
\]
\end{thm}

\end{document}